\documentclass[a4paper,12pt,leqno,twoside]{article}

\usepackage[english]{babel}
\usepackage[latin1]{inputenc}

\usepackage{amsfonts}
\usepackage{amssymb}
\usepackage{amsthm, amsmath}
\usepackage{eucal}
\usepackage{mathrsfs}
\allowdisplaybreaks[0]
\usepackage[hypertexnames=false]{hyperref}

\newcommand\testopari{\sc\small Luca Calatroni and Pierluigi Colli}
\newcommand\testodispari{\sc\small Allen-Cahn Eq. with Singular Potentials and Dynamic Bound. Conds.}
\usepackage[usenames,dvipsnames]{color}
\def\rosso #1{{\color{red}#1}}
\def\blu #1{{\color{blue}#1}}
\def\cmag #1{{\color{magenta}#1}}
\let\rosso\relax
\let\blu\relax
\let\cmag\relax

\newcommand{\R}{\mathbb{R}}

\newenvironment{pr}{\begin{proof}[\textbf{Proof}]}{\end{proof}}
\newtheorem{teor}{Theorem}[section]
\newtheorem{lemm}[teor]{Lemma}

\newtheorem{propos}[teor]{Proposition}

\newtheorem{obs}[teor]{Remark}

\providecommand{\norm}[1]{\left\lVert#1\right\rVert}

\begin{document}
\renewcommand{\baselinestretch}{1.05}
\thispagestyle{empty}

                   \begin{center}
{\large \bf GLOBAL SOLUTION TO THE ALLEN-CAHN EQUATION\\[0.1cm]
WITH SINGULAR POTENTIALS\\[0.3cm]
AND DYNAMIC BOUNDARY CONDITIONS}\\[1cm]
                   \end{center}

                   \begin{center}
                  {\sc Luca Calatroni}\\
 Cambridge Centre for Analysis, University of Cambridge\\
 Wilberforce Road, CB3 0WA, Cambridge, United Kingdom \\
 (lc524@cam.ac.uk) \\[0.4cm]
{\sc Pierluigi Colli}\\
 Dipartimento di Matematica ``F. Casorati'', Universit\`a di Pavia\\
 Via Ferrata 1, 27100, Pavia, Italy \\
 (pierluigi.colli@unipv.it)
       \end{center}
\vskip1cm
{\bf Abstract.} We prove well-posedness results for the solution to an initial and boundary-value problem for an Allen-Cahn type equation describing the phenomenon of phase transitions
for a material contained in a bounded and regular domain. The \emph{dynamic} boundary conditions for the order parameter have been recently proposed by some physicists to account for interactions with the walls. We show our results using suitable regularizations of the nonlinearities of the problem and performing some \emph{a priori} estimates which allow us to pass to the limit thanks to compactness and monotonicity arguments.\\[.2cm]
\centerline{--------------------------------------------------------------}
\\[.2cm]
{
}\\[.2cm]
\centerline{--------------------------------------------------------------}
\\[.5cm]

\section{Introduction}\label{sec:int}

The Allen-Cahn equation was originally introduced \blu{in \cite{al}} as a phenomenological model for anti-phase domain coarsening in a binary alloy. It has been subsequently applied to a wide range of other different problems such as the motion by mean curvature flows (cf., e.g., \cite{feng}) and the crystal growth (cf., e.g., \cite{whe}). In this paper we deal with a physical model which exploits this equation to describe the phenomenon of phase transitions (cf. \cite{brok}).

Let us consider a material occupying at any time $t\geq 0$ a bounded and connected domain $\Omega\subset \R^d,\ d=2,3$ with a smooth boundary $\Gamma:=\partial\Omega$. We suppose  that the material can exhibit two different phases. The following semilinear parabolic partial differential equation proposed by Allen and Cahn in \cite{al} models the time evolution of the order parameter $u$ in the \emph{isothermal} case:
\begin{equation} \label{allencahn}
\partial_t u - \Delta u+W'(u)=0\quad\text{ in }\Omega\times(0,+\infty).
\end{equation}
The nonlinearity $W$ is the potential associated with the phase configuration. In the literature (see, e.g., \cite{el, gal, schimp}) $W$ is frequently assumed to be a double-well potential, i.e. $W(u)=\frac{1}{4}(u^2-1)^2$. Usually, the values $u=\pm 1$ represent the pure states of the material, whereas the values inside the open interval $(-1,1)$ describe a local presence of a mixture of the two states.

In this work we prove well-posedness results for a problem similar to
\eqref{allencahn}, where more general nonlinearities appear. In fact, in our
analysis the potential $W$ may be non-smooth and its domain \rosso{can} be
different from the whole real line\rosso{: thus, we admit the so-called}
\emph{singular} potentials. Physically meaningful examples of these
potentials \rosso{prescribe} $W$ only inside the interval $[-1,1]$, while \rosso{for} the values outside the interval \rosso{one has $W(u)=+\infty$.} Here, we endow our problem with suitable boundary conditions and we
deal with two singular potentials, the former applied to $u$ \blu{in} the interior
of $\Omega$, the latter acting on the trace defined on
$\Gamma$, linking them by a suitable compatibility condition. Each potential will be the sum of a principal part, assumed to be convex and with possible barriers (like the values $\pm 1$), and of a non-convex perturbation. Moreover, the derivative $W'$ in \eqref{allencahn} has to be understood as a subdifferential $\partial W$ and then the \blu{equation} has to be read as a differential inclusion, due to the possibly multivalued character of $\partial W$.

A well-known generalization to the Allen-Cahn equation, which models the evolution of $u$ taking into account also the possible variations of the relative temperature of the system (\emph{non-isothermal} case), has been proposed by Caginalp (see \cite{cag}) for melting-solidification processes in several classes of materials. It is obtained by coupling an equation like \eqref{allencahn} for the order parameter with a heat balance equation for the temperature.

The mathematical literature regarding the Allen-Cahn equation and the Caginalp system is rather vast. In particular, well-posedness results can be found in \cite{schimp}, while the analysis of the dissipative dynamical system generated by these equations has been carried out in a large number of papers. In \cite{gras,gras2} singular potentials are considered  and theorems about the existence of global and/or exponential attractor are proved. We refer the reader also to \cite{cahn} and \cite{miran}, where such models are linked to
the Cahn-Hilliard equation \blu{describing the \emph{phase-separation} phenomenon}, i.e. the progress toward thermodynamic equilibrium of an initial mixture of the phases of the material.

The Allen-Cahn equation and the Caginalp system are usually coupled with homogeneous Neumann boundary conditions for the order parameter. Such conditions are meant to represent the orthogonality of the interface to the boundary and the absence of mass flux. Recently, physicists have introduced 
the so-called \emph{dynamic boundary conditions}, using this terminology to underline the fact that the kinetics of the process, i.e. the term $\partial_t u$, appears explicitly in the boundary conditions. This particular choice has the advantage of taking into account the interaction of the components of the system with the walls (i.e., within $\Gamma$). Dynamic 
conditions can be associated with equation \eqref{allencahn} arguing as described in \cite{gal}. The main idea is to consider a free energy functional \blu{in} the interior of $\Omega$ and another one defined only on $\Gamma$ in which a sufficiently smooth potential $W_b$ appears. One writes 
that the density on the boundary of the total free energy $E$ relaxes towards 
\blu{the} equilibrium with a rate \blu{proportional} to the Fr\'echet-derivative of $E$ with respect to $u$; then, it is possible to infer that the 
system is subject to the following dynamic boundary condition \begin{equation} \label{dynbc}
\frac{1}{\delta}\partial_t u=\alpha\Delta_\Gamma u-\partial_n u-\gamma u-W'_b (u)\quad\text{on }\Gamma\times (0,+\infty),
\end{equation}
\noindent where the parameters $\delta,\, \alpha,\, \gamma$ are positive, $\Delta_\Gamma$ is the Laplace-Beltrami operator and $\partial_n u$ denotes the outer normal derivative of $u$ on $\Gamma$.

Allen-Cahn equation and Caginalp system endowed with boundary conditions \eqref{dynbc} have already been examined in some papers by assuming that the 
potential $W$ has a polynomially controlled growth at most of degree six. 
Well-posedness results as well as the convergence of the solution to a steady state have been discussed in \cite{chil} by means of the {\L}ojasiewicz-Simon 
technique. In \cite{gal}-\cite{gal2} the authors prove global existence and uniqueness \blu{of the solutions} and analyze \blu{their asymptotic behavior}; they show the existence of a global attractor, as well as of an 
exponential attractor, and convergence to equilibrium as time tends to infi\-nity. In the work of Sprekels and Wu \cite{sprek} some analogous results are 
proved together with a suitable {\L}ojasiewicz-Simon type inequality, under the assumptions that the nonlinear terms are real analytic and satisfy some 
particular growth conditions.

In this paper we are interested to the study of the following problem, where the potentials have been generalized:
\begin{align}
 &\partial_t u-\Delta u+\beta(u)+\pi(u)\ni f\quad \text{in }\Omega\times (0,T) ,\label{eqint}\\
 &v:=u|_{\Gamma},\quad\partial_n u+\partial_t v-\nu\Delta_{\Gamma}v+\beta_{\Gamma}(v)+\pi_{\Gamma}(v)\ni f_{\Gamma}\quad\text{in }\Gamma\times (0,T) ,\label{eqbound}\\
 &u(0)=u_0\, \text{ in }\Omega\quad\text{and }\ v(0)=v_0\, \text{ on }\Gamma .\label{initialcond}
\end{align}
\noindent Here, $T$ is an arbitrary positive and finite final time, $\beta$ and $\beta_{\Gamma}$ are maximal monotone and \rosso{possibly} non-smooth graphs, while $\pi$ and $\pi_{\Gamma}$ are non-monotone smooth perturbations. Moreover, we have introduced the nonnegative \blu{parameter} $\nu$ in \eqref{eqbound} and \blu{added} the forcing terms $f$ and $f_{\Gamma}$. We consider nonlinear terms with the following regularities: $\pi$ and $\pi_{\Gamma}$ are Lipschitz continuous functions, while the graphs $\beta$ and $\beta_{\Gamma}$ are naturally subdifferentials of convex functions. Therefore, the equations \eqref{eqint} and \eqref{eqbound} have to be read as differential inclusions. Recalling equation \eqref{allencahn} endowed with the boundary conditions \eqref{dynbc}, we observe that our choice corresponds to split $W'$ and $W'_b$ as $W'=\partial W=\beta+\pi$ and $W'_b=\partial W_b =\beta_\Gamma+\pi_\Gamma$. A relevant example for $\beta$ and $\beta_\Gamma$ is provided by  $\beta(u)=\partial I_{[-1,1]}(u)$ and $\beta_\Gamma (v)=\partial I_{[-1,1]}(v)$, where the function $I_{[-1,1]}:\R\rightarrow (-\infty,+\infty]$ is the indicator function of the interval $[-1,1]$. We want to obtain well-posedness results for the problem \eqref{eqint}-\eqref{initialcond} on any finite time interval $(0, T)$: our arguments are inspired by those employed by Gilardi, Miranville and Schimperna in \cite{gil} when dealing with the Cahn-Hilliard equation.

\subsection{Outline of the paper}
We start our investigation by introducing an equivalent abstract formulation of \eqref{eqint} which takes into account both the boundary condition
\eqref{eqbound} and the initial conditions \eqref{initialcond}. After defining the framework in which we will develop our work, we state the
continuous dependence on data result (Theorem \ref{contdep}) and, as a
corollary, the uniqueness result (Theorem \ref{uniqueness}), both holding under quite general assumptions on the data. Existence of a weak solution is
proved under more restrictive, though natural, assumptions (Theorem \ref{exist}). In particular,  we assume that the main nonlinearities $\beta$
and $\beta_\Gamma$ satisfy a proper compatibility condition. The existence proof is given using a two-level approximating scheme. Our argument is
\rosso{detailed} in Section \ref{exis:te} and runs as follows. First, as the nonlinearities $\beta$ and $\beta_\Gamma$ of our problem could be
multivalued, we regularize them in order to deal with single-valued and more regular functions. This is what we do at the first level of our
approximation \blu{by using Yosida approximations depending on a small parameter $\varepsilon$}. Then, \rosso{with the help of} a suitable finite-dimensional approximating scheme and thanks to some
technical results about the regularizations, we are able to find some a priori estimates on our $\varepsilon$-approximation. By such estimates we can
pass to the limit as $\varepsilon\searrow 0$ taking the maximal monotonicity of the graphs into account and by virtue of compactness and monotonicity
results. Lastly, Section \ref{case nu0} is dedicated  to the analysis of the problem in the special case \blu{where} the parameter $\nu$ is equal to $0$\rosso{: here, besides the well-posedness issues, we can prove a convergence result for the solutions of the problem with $\nu >0 $ as $\nu \searrow 0$.}

\subsection{Notation}
 We conclude this section stating a general rule regarding the notation in the paper. In order to avoid a boring notation, throughout the paper the
symbols $c$ and $C$ stand for different constants which depend generally on the domain and on the norms of the functions involved. In particular, such
constants are independent \rosso{of the coefficient} $\nu$ and \rosso{the approximation parameter}
$\varepsilon$ we introduce in next section. On the other hand, a notation
like $c_\delta$ allows the constant to depend in addition on \rosso{some} parameter $\delta$. Hence, the meaning of $c$, $C$ and $c_\delta$ might change from
line to line and even in the same chain of inequalities.

\section{Main results} \label{main:res}
\setcounter{equation}{0}

Let $\Omega$ be a bounded domain in $\R^d$, $d>1$, with smooth boundary $\Gamma:=\partial\Omega$. We denote by $\partial_n$ the outward normal derivative on $\Gamma$. Let $|\Omega|$ stand for the Lebesgue measure of $\Omega$ and let $|\Gamma|$ be the $(d-1)$-dimensional measure of $\Gamma$. Given a final time $T$, we set for convenience
\begin{align}
 Q_t:=\Omega\times(0,t)\quad&\text{and}\quad \Sigma_t:=\Gamma\times(0,t)\ \text{for every }t\in (0,T], \notag\\
 Q:=Q_T\quad&\text{and}\quad\Sigma:=\Sigma_T. \notag
\end{align}
We introduce the functions $\hat{\beta}, \hat{\beta_{\Gamma}},\pi,\pi_{\Gamma}$ and the positive parameter $\nu$ (the case $\nu$=0 will be treated separately) satisfying the following conditions:
\begin{align}
& \hat{\beta},\hat{\beta_{\Gamma}}:\R\rightarrow[0,+\infty)\ \text{are convex, proper, l.s.c., and}\ \hat{\beta}(0)=\hat{\beta_{\Gamma}}(0)=0, \label{convfunct}\\
& \pi ,\, \ \pi_{\Gamma}\ \text{are Lipschitz continuous with Lipschitz constants } L, \, L_\Gamma\text{ respectively}, \label{lipschitz}\\
& \nu> 0. \label{nonnegatpar}
\end{align}
We can therefore consider the subdifferentials $\beta$ and $\beta_{\Gamma}$ in $\R\times\R$ of the functions $\hat{\beta}, \hat{\beta_{\Gamma}}$:
\begin{equation} \label{subdiff}
\beta:=\partial\hat{\beta}\quad\text{and}\quad\beta_{\Gamma}:=\partial\hat{\beta_{\Gamma}},
\end{equation}
\noindent which turn out to be maximal monotone graphs in $\R\times\R$ (we refer to \cite{bar} and \cite{bre} for basic definitions and properties of maximal monotone operators). We note that \eqref{convfunct} entails $\beta(0)\ni 0$ and $\beta_{\Gamma}(0)\ni 0$. In order to let these graphs be as general as possible, we allow them to have \emph{effective domains}, respectively denoted by
\begin{equation*}
D(\beta)=\left\{r\in\R: \beta(r)\neq\emptyset\right\}\quad\text{and}\quad D(\beta_{\Gamma})=\left\{r\in\R:\beta_\Gamma(r)\neq\emptyset\right\},
\end{equation*}
\rosso{possibly} different from the whole real line. \blu{For $\varepsilon>0$} we introduce also the following notation for any maximal monotone graph $\alpha:\R\rightarrow2^{\R}$ (see, e.g., \cite[p.~28]{bre})
\begin{align}
& \alpha^{\circ}(r)\ \text{is the element of }\alpha(r)\ \text{having minimum modulus}, \label{yosida}\\
& \alpha^{Y}_{\varepsilon}:=\varepsilon^{-1}(I-(I+\varepsilon\alpha)^{-1})\ \text{is the Yosida \blu{approximation} of}\ \alpha\text{ of parameter }\varepsilon. \notag
\end{align}
Just for convenience, we still use the symbol $\alpha$ (and, e.g., $\alpha^Y_\varepsilon$ as a particular case) to represent the induced maximal monotone operator in  $L^2(\Omega)$ (more precisely, from $L^2(\Omega)$ to $L^2(\Omega)$), doing the same also for the operators induced in $L^2(\Gamma), L^2(Q)$ and $L^2(\Sigma)$.

As the graphs $\beta$ and $\beta_{\Gamma}$ might be multivalued, in our formulation we introduce two selections $\xi$ and $\xi_{\Gamma}$ from $\beta$ and $\beta_{\Gamma}$, respectively. Therefore, the problem \eqref{eqint}-\eqref{initialcond} can be stated as follows:
\begin{align}
& \partial_t u-\Delta u+\xi+\pi(u)=f\quad \text{a.e. in}\ Q  ,\label{analform}\\
& v:=u|_{\Gamma},\quad\partial_n u+\partial_t v-\nu\Delta_{\Gamma}v+\xi_\Gamma+\pi_{\Gamma}(v)= f_{\Gamma}\quad \text{a.e. on}\ \Sigma  ,\label{boundan}\\
& \xi\in\beta(u)\ \text{a.e. in }Q\quad\text{and}\quad\xi_{\Gamma}\in\beta_{\Gamma}(v)\ \text{a.e. on }\Sigma , \\
& u(0)=u_0\text{ a.e. in }\Omega\quad\text{and }\ v(0)=v_0\text{ a.e. on }\Gamma. \label{initcondprob}
\end{align}

We point out the functional spaces we will use \blu{in the following} and set
\begin{align}
& V:=H^1(\Omega),\quad H:=L^2(\Omega),  \label{VH} \\
&\hskip1.5cm V_{\Gamma}:=H^1(\Gamma),\quad H_{\Gamma}:=L^2(\Gamma),\quad \mathcal{V}:=\left\{z\in V:z|_{\Gamma}\in V_{\Gamma}\right\}. \notag
\end{align}
Note that the inclusion $V\subset H$ is compact and dense. Then we identify $H$ with its dual space $H^{*}$,
so that $V\subset H\subset V^{*}$ with dense and compact injections, and therefore $(V,H,V^{*})$ is the standard Hilbert triplet. The same considerations hold for $V_\Gamma$ and $H_\Gamma$. Let $(\cdot,\cdot)$ and $\norm{\cdot}_H$ be the inner product and the induced norm in $H$. We also denote by $\norm{\cdot}_V, \norm{\cdot}_{H_\Gamma}, \norm{\cdot}_{V_\Gamma}$ the norms in the corresponding spaces. The symbol $\langle\cdot,\cdot\rangle$ stands for the duality pairing between $V^{*}$ and $V$. We endow the space $\mathcal{V}\subset V$ with the following inner product and the induced graph norm:
\begin{equation}
 (w,z)_{\mathcal{V}}:=(w,z)_V+(w|_\Gamma,z|_\Gamma)_{V_\Gamma} , \quad
 \norm{w}^2_{\mathcal{V}}=\norm{w}^2_V+\norm{w|_{\Gamma}}^2_{V_\Gamma}, \quad w,z\in\mathcal{V}   . \label{normvcors}
\end{equation}

We present now the variational formulation of the problem \eqref{analform}-\eqref{initcondprob} in the abstract setting introduced above. In order to obtain it, we formally multiply the equation \eqref{analform} by test functions and integrate by parts on $\Omega$ taking the
boundary conditions \eqref{boundan} into account. We have to make precise the regularity of all the components of the problem: the results of
continuous dependence and uniqueness of the solution hold under quite general
assumptions, while further requirements have to be specified in order to
prove \rosso{the existence result}.
So, just to start, we ask the data to fulfill the conditions
\begin{equation}       \label{datacond}
f\in L^2(0,T;H),\quad f_{\Gamma}\in L^2(0,T;H_{\Gamma}),\quad u_0\in H\quad\text{and}\quad v_0\in H_\Gamma
\end{equation}
\noindent and look for a quadruplet $(u,v,\xi,\xi_{\Gamma})$ such that
\begin{align}
& u\in L^2(0,T;V)\cap H^1(0,T;V^{*}), \label{reguniqu}\\
& v:=u|_{\Gamma}\in L^2(0,T;V_{\Gamma})\cap H^1(0,T;V^{*}_{\Gamma}), \label{reguniqv} \\
& \xi\in L^2(0,T;H)\rosso{{}\equiv L^2(Q)} \ \text{ and }\ \xi\in\beta(u)\text{ a.e. in }Q ,\label{reguniqcsi}\\
& \xi_{\Gamma}\in L^2(0,T;H_{\Gamma})\rosso{{}\equiv L^2(\Sigma)}\ \text{ and }\ \xi_\Gamma\in\beta_\Gamma(v)\text{ a.e. on }\Sigma,\label{reguniqcsigamma}\\
& u(0)=u_0,\ v(0)=v_0  \label{dataconduniq}
\end{align}
\noindent and satisfying for a.a. $t\in (0,T)$ and every $z\in\mathcal{V}$
\begin{align}
& \langle\partial_t u(t),z\rangle+\int_{\Omega}\nabla u(t)\cdot\nabla z +\int_{\Omega}(\xi(t)+\pi(u(t)))z+ {_{V^{*}_{\Gamma}}}\langle\partial_t v(t),z\rangle_{V_{\Gamma}} \label{varform}\\
& +\nu\int_{\Gamma}\nabla_{\Gamma}v(t)\cdot\nabla_{\Gamma}z+\int_{\Gamma}(\xi_{\Gamma}(t)+\pi_{\Gamma}(v(t)))z=\int_{\Omega}f(t)z+\int_{\Gamma}f_{\Gamma}(t)z .\notag
\end{align}

\begin{obs} \label{nottrace}
With a little abuse of notation we will use the same letter $z$ to denote the function $z\in\mathcal{V}$ and its trace $z|_\Gamma \in V_{\Gamma} $. If $z\in L^2(0,T;V)$, the meaning of $\partial_t z$ on the boundary has to be understood as $\partial_t (z|_{\Gamma})$ which exists at least in a distributional sense.
\end{obs}

\begin{obs}
We observe that regularities \eqref{reguniqu} and \eqref{reguniqv} imply that
\begin{equation*}
u\in C^0([0,T];H)\quad\text{and}\quad v\in C^0([0,T];H_\Gamma).
\end{equation*}
Thus, in particular, the values $u(t)$ and $v(t)$ make sense in $H$ and $H_\Gamma$, respectively, for all $t\in [0,T]$, and this gives a meaning to the initial conditions \eqref{dataconduniq} in such spaces.
\end{obs}

In the framework specified by \eqref{datacond}-\eqref{reguniqcsigamma}, the following continuous dependence result holds true, together with a uniqueness one as its immediate consequence.
\begin{teor} [Continuous dependence on the data] \label{contdep}
Assume \eqref{convfunct}-\eqref{subdiff},
\eqref{datacond} and let $(u_1,v_1)$,
$(u_2,v_2)$ be the first two components
of two different solutions of the problem
\eqref{reguniqu}-\eqref{varform}, each one
corresponding to a set of data $u_{0,i},v_{0,i},f_i,f_{\Gamma,i}$ with $i=1,2$. Then, there is a constant $c$ such that
\begin{align}\label{contdip}
& \norm{u_1-u_2}^2_{L^{\infty}(0,T;H)}+\norm{\nabla(u_1-u_2)}^2_{L^2(Q)}\\
&+ \norm{v_1-v_2}^2_{L^{\infty}(0,T;H_{\Gamma})}
+\nu\norm{\nabla_{\Gamma}(v_1-v_2)}^2_{L^2(\Sigma)}\notag\\
& \leq c\, \Big\{\norm{u_{0,1}-u_{0,2}}^2_H
+ \norm{v_{0,1}-v_{0,2}}^2_{H_\Gamma}\notag\\
&\hskip1.5cm {}+\norm{f_1-f_2}^2_{L^2(0,T;H)}+\norm{f_{\Gamma,1}-f_{\Gamma,2}}^2_{L^2(0,T;H_{\Gamma})}\Big\}. \notag
\end{align}

\end{teor}

\begin{teor} [Uniqueness] \label{uniqueness}
Under the same assumptions as in Theorem \ref{contdep}, then any two solutions of the problem \eqref{reguniqu}-\eqref{varform} coincide.
\end{teor}

In order to prove the existence of the solution, we have to enforce our assumptions requiring:
\begin{align}
& u_0\in V,\quad v_0=u_0|_{\Gamma}\in V_{\Gamma} \label{initcondex}\\
& \hat{\beta}(u_0)\in L^1(\Omega)\ \ \text{and}\ \ \hat{\beta}_{\Gamma}(v_0)\in L^1(\Gamma). \label{betacondex}
\end{align}
Of course, \eqref{initcondex} is equivalent to $u_0\in\mathcal{V}$ and the trace of $u_0$ on the boundary coincides with the initial value $v_0$. Moreover, we introduce also a compatibility condition on the main nonlinearities $\beta$ and $\beta_{\Gamma}$. Namely, we assume that:
\begin{equation} \label{domcont}
D(\beta)\supseteq D(\beta_{\Gamma})
\end{equation}
\noindent and that two real constants $\eta$ and $C_{\Gamma}$ exist such that
\begin{align}
& \eta>0,\quad C_{\Gamma}\geq 0, \quad 
 |\beta^{\circ}(r)|\leq\eta|\beta^{\circ}_{\Gamma}(r)|+C_{\Gamma}\ \text{  for every  }r\in D(\beta_{\Gamma}). \label{condcomp}
\end{align}

We state now our main result, whose proof is given in Section~\ref{exis:te}.

\begin{teor} [Existence]  \label{exist}
Assume \eqref{convfunct}-\eqref{subdiff}, \eqref{datacond} and \eqref{initcondex}-\eqref{condcomp}. Then, there exists a quadruplet $(u,v,\xi,\xi_{\Gamma})$ solving problem \eqref{reguniqu}-\eqref{varform}. Moreover, $u$ and $v$ fulfill
\begin{equation}
 u\in L^{\infty}(0,T;V)\cap H^1(0,T;H) , \ \quad
 v\in L^{\infty}(0,T;V_\Gamma)\cap H^1(0,T;H_\Gamma).\notag
\end{equation}
\end{teor}

In Section \ref{case nu0} we prove that very similar continuous dependence and existence results hold in the case $\nu=0$ \rosso{(see Theorems~\ref{uniqnu0} and~\ref{exist0})}. As we will point out, in order to do this we have to change our choice of the functional spaces on the boundary and extend in a suitable \rosso{way} the operator related to $\beta_\Gamma$. Moreover, we will be able to show that the solution of the problem for $\nu>0$ converges in a suitable topology to the solution of the problem with $\nu=0$.

\section{Continuous dependence on the data} \label{con:dep}
\setcounter{equation}{0}

This section contains the proofs of Theorems~\ref{contdep} and~\ref{uniqueness}. In order to prove Theorem~\ref{contdep} we consider two
different solutions of the problem \eqref{reguniqu}-\eqref{varform} coming each from a different set of data and we label their components with
subscripts 1 and 2. We write \eqref{varform} for both solutions and take the difference, defining $u$ as $u:= u_1-u_2$ and doing the same for $v:=v_1-v_2,
$ $\xi:=\xi_1-\xi_2,$ $\xi_{\Gamma}:=\rosso{\xi_{\Gamma, 1}-\xi_{\Gamma , 2}}$. Then, we write such a difference at time $t=s$, we choose as test-function $z\in\mathcal{V}$ the function $u(s)=u_1 (s)-u_2 (s)$ itself and integrate what we get over $(0,t)$, where $t\in(0,T]$ is arbitrary. We obtain
\begin{align}
& \frac{1}{2}\norm{u(t)}^2_H + \frac{1}{2}\norm{v(t)}^2_{H_\Gamma} -\frac{1}{2}\norm{u_{0,1}-u_{0,2}}^2_H -\frac{1}{2}\norm{v_{0,1}-v_{0,2}}^2_{H_\Gamma} \label{contdip1}\\
&\quad {}+\int_{Q_t}|\nabla u|^2+\nu\int_{\Sigma_t}|\nabla_{\Gamma}v|^2 +\int_{Q_t}\xi u+\int_{\Sigma_t}\xi_{\Gamma}v \notag\\
& =\int_{Q_t}(\pi(u_2)-\pi(u_1))u+ \int_{\Sigma_t}(\pi_{\Gamma}(v_2)-\pi_{\Gamma}(v_1))v \notag\\
&\quad {}+\int_{Q_t}(f_1-f_2)u +\int_{\Sigma_t}(f_{\Gamma,1}-f_{\Gamma,2})v. \notag
\end{align}
The last two terms on the left-hand side of \eqref{contdip1} are nonnegative, due to the monotonicity properties of $\beta$ and $\beta_\Gamma$. As far as the right-hand side is concerned, we can estimate the first two terms \rosso{with the help of} \eqref{lipschitz}:
\begin{equation} \label{estlip}
 \int_{Q_t}(\pi(u_2)-\pi(u_1))u+\int_{\Sigma_t}(\pi_{\Gamma}(v_2)-\pi_{\Gamma}(v_1))v\leq L\int_{Q_t}|u|^2+L_{\Gamma}\int_{\Sigma_t}|v|^2 \rosso{.}
\end{equation}
\noindent \rosso{For the other two terms, we simply use} Young's inequality:
\begin{align}
\label{contdip2}
&\int_{Q_t}(f_1-f_2)u +\int_{\Sigma_t}(f_{\Gamma,1}-f_{\Gamma,2})v \\
&\leq\frac{1}{2}\int_{Q_t}|f_1-f_2|^2+\frac{1}{2}\int_{Q_t}|u|^2 +\frac{1}{2}\int_{\Sigma_t}|f_{\Gamma,1}-f_{\Gamma,2}|^2+\frac{1}{2}\int_{\Sigma_t}|v|^2. \notag
\end{align}
Thus, \rosso{thanks to} \eqref{contdip1}-\eqref{contdip2} we have
\begin{align}
& \norm{u(t)}^2_H+\norm{v(t)}^2_{H_\Gamma}+\norm{\nabla u}^2_{L^2(Q_t)}+\nu\norm{\nabla_\Gamma v}^2_{L^2(\Sigma_t)}\notag \\
& \leq \int_{Q_t}|f_1-f_2|^2 +\int_{\Sigma_t}|f_{\Gamma,1}-f_{\Gamma,2}|^2 +(2L+1)\int_{Q_t}|u|^2+(2L_\Gamma+1)\int_{\Sigma_t}|v|^2 \notag \\
& +\norm{u_{0,1}-u_{0,2}}^2_H+\norm{v_{0,1}-v_{0,2}}^2_{H_\Gamma}. \notag
\end{align}
We can now apply the Gronwall Lemma and then take the supremum over the interval $(0,T)$ thus finding $L^\infty$-in time estimates for the left-hand side of this inequality. Scaling the constants, we finally get the estimate \eqref{contdip}.

For the proof of Theorem~\ref{uniqueness}, we note that Theorem~\ref{contdep} yields the uniqueness of the first two components of the solution, that is, $u_1 =u_2$ and  $v_1 =v_2$. Then, with the help of
\eqref{varform} we easily obtain
\begin{equation} \label{pc21}
\int_{\Omega}(\xi_1 - \xi_2) (t)z+ \int_{\Gamma}(\xi_{\Gamma,  1 }- \xi_{\Gamma,  2})
(t)z = 0 \quad \text{for all } \, z \in \mathcal{V}, \  \text{ for a.a. } \,
t \in (0,T).
\end{equation}
Taking first $z\in H^1_0(\Omega) (\subset \mathcal{V}) $, we easily deduce that $\xi_1   = \xi_2$; then, \blu{by rewriting \eqref{pc21} for a general test function $z\in \mathcal{V}$, we can conclude that $\xi_{\Gamma,  1 }= \xi_{\Gamma,  2}$.}

\section{Existence} \label{exis:te}
\setcounter{equation}{0}

This section is devoted to the proof of  Theorem~\ref{exist}. We fix a parameter $\varepsilon\in (0,1)$ and consider an approximating problem depending on $\varepsilon$, \blu{obtained} by regularizing the nonlinearities $\beta$ and $\beta_{\Gamma}$. To this aim, we define the functions $\beta_{\varepsilon}$ and $\beta_{\Gamma,\varepsilon}:\R\rightarrow\R$ by
\begin{align}
& \beta_{\varepsilon}(r):=\beta^Y_{\varepsilon}(r)\ \ \ \text{for }r\in\R , \label{yosidaint}\\
& \beta_{\Gamma,\varepsilon}(r):=(\beta_\Gamma)^Y_{\varepsilon\eta}(r)\ \ \ \text{for }r\in\R .\label{yosidabound}
\end{align}
Here, we have used the notation \eqref{yosida} for the Yosida approximations and the parameter $\eta$ is the fixed constant given by \eqref{condcomp}. The choice of a function $\beta_{\Gamma,\varepsilon}$ as above is justified by the forthcoming Lemma~\ref{growcondeps} which ensures a compatibility condition similar to \eqref{condcomp}.

For convenience\blu{,} we define also the functions $\hat{\beta}_{\varepsilon}, \hat{\beta}_{\Gamma,\varepsilon}:\R\rightarrow\R$ by
\begin{equation} \label{primitbeta}
\hat{\beta}_{\varepsilon}(r):=\int_0^r\beta_{\varepsilon}(s)ds\quad\text{and }\quad
\hat{\beta}_{\Gamma,\varepsilon}(r):=\int_0^r\beta_{\Gamma,\varepsilon}(s)ds\ \ \text{for }r\in\R.
\end{equation}
Let us recall that the Yosida $\lambda$-regularization of a maximal monotone operator is monotone and Lipschitz continuous with $1/\lambda$ as Lipschitz constant (see \cite[Propositions~2.6 and~2.7]{bre}). Thus, such a property holds in particular for $\beta_{\Gamma,\varepsilon}$ with $1/\varepsilon\eta$ as Lipschitz constant. Moreover, both $\beta_{\varepsilon}$ and $\beta_{\Gamma,\varepsilon}$ vanish at 0. It follows that $\hat{\beta}_{\varepsilon}$ and $\hat{\beta}_{\Gamma,\varepsilon}$ are
nonnegative convex functions with (at most) quadratic \blu{growth}. We finally point out a general property of Yosida approximations of every maximal monotone operator (see \cite[Prop.~2.6, p.~28]{bre})\blu{:}
    \begin{align}
|\beta_\varepsilon(r)|\leq|\beta^{\circ}(r)|\ \text{  for   }
\, r\in D(\beta), \quad
\ |\beta_{\Gamma,\varepsilon}(r)|\leq |\beta^{\circ}_\Gamma(r)| \ \text{  for   }\, r\in D(\beta_\Gamma), \label{ineqeps}\\
  0 \leq \hat{\beta}_\varepsilon(r) \leq \hat{\beta}(r) \ \text{  and  }\ 0 \leq \hat{\beta}_{\Gamma,\varepsilon}(r)\leq \hat{\beta_\Gamma}(r) \
\text{  for all  } \, r\in \R.
\notag
     \end{align}
Our approximating problem consists in finding a function $u_{\varepsilon}$ satisfying
\begin{align}
& u_{\varepsilon}\in L^\infty(0,T;V)\cap H^1(0,T;H) \label{regueps}\\
& u_{\varepsilon}|_{\Gamma}\in L^\infty(0,T;V_{\Gamma})\cap H^1(0,T;H_{\Gamma}) \label{reguepsboun}\\
& u_{\varepsilon}(0)=u_0 \label{datacondeps}
\end{align}
\noindent and solving for a.a. $t\in(0,T)$ the variational equality
\begin{align}
& \int_{\Omega}\partial_t u_{\varepsilon}(t)z+\int_{\Omega}\nabla u_{\varepsilon}(t)\cdot\nabla z +\int_{\Omega}(\beta_{\varepsilon}+\pi)(u_{\varepsilon}(t))z+\int_{\Gamma}\partial_t u_{\varepsilon}(t)z \label{varformeps} \\
& +\nu\int_{\Gamma}\nabla_{\Gamma}u_{\varepsilon}(t)\cdot\nabla_{\Gamma}z+\int_{\Gamma}(\beta_{\Gamma,\varepsilon}+\pi_{\Gamma})(u_{\varepsilon}(t))z
=\int_{\Omega}f(t)z+\int_{\Gamma}f_{\Gamma}(t)z \notag
\end{align}
\noindent for every $z\in\mathcal{V}$, where, referring to Remark \ref{nottrace}, from now on we simply write $u_\varepsilon$ instead of $u_\varepsilon|_\Gamma$ \blu{to indicate the trace of $u_\varepsilon$ on the boundary}.

We state our well-posedness result for the above problem.
\begin{propos} \label{wellposedprobleps} Assume \eqref{lipschitz}, \eqref{datacond} and \eqref{initcondex}. Let $\beta_{\varepsilon}$ and $\beta_{\Gamma,\varepsilon}:\R\rightarrow\R$ be the functions introduced in \eqref{yosidaint}-\eqref{yosidabound}.  Then, there exists a unique function $u_{\varepsilon}$ satisfying \eqref{regueps}-\eqref{varformeps}.
\end{propos}

The uniqueness proof follows as a particular case of Theorem~\ref{uniqueness}. In order to show the existence of the solution, we \blu{need to} introduce an approximation of our $\varepsilon$-problem in a finite dimensional space using a standard Faedo-Galerkin technique. Hence, we can show that a more regular discrete solution $u^n_\varepsilon$ exists by virtue of well-known results
for systems of ordinary differential equations. We can thus find the desired $u_\varepsilon$ by passing to the limit as $n\nearrow +\infty$ and exploiting some compactness and monotonicity results (see \cite{gil} and \cite{luca} for further details).

\subsection{Estimates}

We perform now some a priori estimates on $u_\varepsilon, \beta_\varepsilon(u_\varepsilon), \beta_{\Gamma,\varepsilon}(u_\varepsilon)$ in order to solve the \rosso{problem \eqref{reguniqu}-\eqref{varform}}
by letting $\varepsilon$ tend to 0. All such estimates can be rigorously justified at the level of the Galerkin scheme (see \cite{luca}).

\subsubsection{First a priori estimate}
We choose $z= u_{\varepsilon}(t)$ in \eqref{varformeps}: we are allowed to do this thanks to \rosso{the} regularities \eqref{regueps}-\eqref{reguepsboun}. With the help of Young's inequality we obtain
\begin{align}
 \frac{1}{2}\frac{d}{dt}\norm{u_\varepsilon(t)}^2_H+\frac{1}{2}\frac{d}{dt}\norm{u_\varepsilon(t)}^2_{H_\Gamma}
+ \norm{\nabla u_\varepsilon(t)}^2_H
+\nu\norm{\nabla_\Gamma u_\varepsilon(t)}^2_{H_\Gamma} \hskip3.5cm \label{firstestn}\\
+\int_\Omega\beta_\varepsilon(u_\varepsilon(t))u_\varepsilon(t)
+\int_\Gamma\beta_{\Gamma,\varepsilon}(u_\varepsilon(t))u_\varepsilon(t)
\leq -\int_\Omega \pi(u_\varepsilon(t))u_\varepsilon(t)-\int_\Gamma
\pi_\Gamma(u_\varepsilon(t))u_\varepsilon(t)
\notag\\
{}+\frac{1}{2}\norm{f(t)}^2_H+\frac{1}{2}\norm{f_\Gamma(t)}^2_{H_\Gamma}+\frac{1}{2}\norm{u_\varepsilon(t)}^2_H+\frac{1}{2}\norm{u_\varepsilon(t)}^2_{H_\Gamma} .  \notag
\end{align}
The two terms containing $\beta_\varepsilon$ and $\beta_{\Gamma,\varepsilon}$
are nonnegative since both $\beta_\varepsilon$ and $\beta_{\Gamma,\varepsilon}$
are monotone and null in $0$, so we can avoid considering them. On the other hand, due to \eqref{lipschitz} we have that
\begin{align*}
&-\int_\Omega \pi(u_\varepsilon(t))u_\varepsilon(t) - \int_\Gamma \pi_\Gamma(u_\varepsilon(t))u_\varepsilon(t)\\
& \leq (L+1)\int_\Omega |u_\varepsilon(t)|^2+(L_\Gamma+1)\int_\Gamma |u_\varepsilon(t)|^2+ C,
\end{align*}
where the constant $C$ depends \blu{just} on $|\Omega|, |\Gamma|$ and on the values taken in $0$ by the functions $\pi$ and $\pi_\Gamma$. Then, integrating both sides of \eqref{firstestn} with respect to time and applying the Gronwall lemma,
we easily find that
\begin{align} \label{gronw}
&\norm{u_\varepsilon(t)}^2_H+\norm{u_\varepsilon(t)}^2_{H_\Gamma} + \norm{\nabla u_\varepsilon}^2_{L^2(0,t;H)} +\nu\norm{\nabla_\Gamma u_\varepsilon}^2_{L^2(0,t;H_\Gamma)} \\
&\leq C \left( 1+  \norm{u_0}^2_H + \norm{u_0}^2_{H_\Gamma} +
\norm{f}^2_{L^2(0,T;H)} + \norm{f_\Gamma}^2_{L^2(0,T;H_\Gamma)} \right) \notag
\end{align}
\noindent for all $t\in (0,T)$. Therefore, in view of \eqref{datacond},
it turns out that
\begin{equation} \label{firstaprepsn}
\norm{u_\varepsilon}_{L^\infty(0,T;H)\cap L^2(0,T;V)}
+ \norm{u_\varepsilon}_{L^\infty(0,T;H_\Gamma) }
+ \nu^{1/2} \norm{u_\varepsilon}_{L^2(0,T;V_\Gamma)}\leq C.
\end{equation}

\subsubsection{Second a priori estimate}
We have now to proceed formally by testing equation \eqref{varformeps} by $z=\partial_t u_{\varepsilon}(t)$ and integrating with respect to $t$. Owing to \eqref{primitbeta}, we obtain
\begin{align}
& \int_{Q_t}|\partial_t u_{\varepsilon}|^2
+ \int_{\Sigma_t}|\partial_t u_{\varepsilon}|^2
+ \frac{1}{2}\int_{\Omega}|\nabla u_\varepsilon(t)|^2
\label{aprioriestim1}\\
&+ \frac{\nu}{2}\int_{\Gamma}|\nabla_{\Gamma} u_\varepsilon (t) |^2
+ \int_{\Omega}\hat{\beta}_\varepsilon(u_\varepsilon(t))
+ \int_{\Gamma}\hat{\beta}_{\Gamma,\varepsilon} (u_\varepsilon(t))
\notag\\
&\leq \frac{1}{2}\int_{\Omega}|\nabla u_0|^2
+ \frac{\nu}{2}\int_{\Gamma}|\nabla_{\Gamma} u_0 |^2
+ \int_{\Omega}\hat{\beta}_\varepsilon(u_0)
+ \int_{\Gamma}\hat{\beta}_{\Gamma,\varepsilon} (u_0)\notag\\
&\quad  + \int_{Q_t} (f -\pi(u_\varepsilon) ) \partial_t u_{\varepsilon}
+ \int_{\Sigma_t} (f_{\Gamma} - \pi_{\Gamma}(u_\varepsilon))\partial_t u_{\varepsilon} . \notag
\end{align}
As all the integrals on the \rosso{left-hand side are nonnegative, we consider the ones on} \blu{the}
right-hand side. Owing to \eqref{datacond}, \eqref{lipschitz} and
\eqref{firstaprepsn}, by applying \blu{Young's} inequality we deduce that
\begin{equation*}
\int_{Q_t} (f -\pi(u_\varepsilon) ) \partial_t u_{\varepsilon}
+ \int_{\Sigma_t} (f_{\Gamma} - \pi_{\Gamma}(u_\varepsilon))\partial_t
u_{\varepsilon} \leq C + \frac{1}{2}\int_{Q_t}|\partial_t u _{\varepsilon}|^2 +
\frac{1}{2}\int_{\Sigma_t}|\partial_t u _\varepsilon|^2 .
\end{equation*}
Next, we examine the terms containing the initial conditions $u _0$.
Thanks to the general property pointed out in \eqref{ineqeps}, we have
\begin{equation*}
\int_{\Omega}\hat{\beta}_\varepsilon(u_0)+\int_{\Gamma}\hat{\beta}_{\Gamma,\varepsilon}(u_0)\leq
\int_{\Omega}\hat{\beta}(u_0)+\int_{\Gamma}\hat{\beta}_{\Gamma}(u_0)
\end{equation*}
and \eqref{initcondex}-\eqref{betacondex} allow us to infer
$$
\frac{1}{2}\int_{\Omega}|\nabla u_0|^2
+ \frac{\nu}{2}\int_{\Gamma}|\nabla_{\Gamma} u_0 |^2 +
\int_{\Omega}\hat{\beta}_\varepsilon(u_0)+\int_{\Gamma}\hat{\beta}_{\Gamma,\varepsilon}(u_0) \leq C.
$$
Then, from \eqref{aprioriestim1} we conclude that
\begin{equation} \label{unibound}
\norm{u_\varepsilon}_{L^{\infty}(0,T;V) \cap H^1(0,T;H)} +\norm{u_\varepsilon}_{H^1(0,T;H_\Gamma)}
+ \nu^{1/2} \norm{u_\varepsilon}_{L^{\infty}(0,T;V_\Gamma)}
\leq C.
\end{equation}

\subsubsection{Preliminary results to the following estimates}
We highlight now some general properties of \rosso{the} Yosida approximation \rosso{that will be} useful to prove \rosso{the subsequent} Lemma~\ref{growcondeps}.
\begin{lemm} \label{yosidacost}
Let $\gamma:\R\rightarrow 2^{\R}$ be a maximal monotone graph and let  $\varepsilon,b>0$. Then
\begin{equation}\label{pr0}
(b\gamma)^Y_\varepsilon=b(\gamma)^Y_{b\varepsilon}.
\end{equation}
\end{lemm}
\begin{pr}
Recalling the definition of Yosida approximation given by \eqref{yosida}, the
equalities
$$
b(\gamma)^Y_{b\varepsilon}=b\frac{1}{b\varepsilon}(I-(I+b\varepsilon\gamma)^{-1})=\frac{1}{\varepsilon}(I-(I+\varepsilon(b\gamma))^{-1})=(b\gamma)^Y_\varepsilon .
$$
enable us to deduce \eqref{pr0}.\end{pr}

\begin{lemm} \label{yosidasum}
Let $\gamma:\R\rightarrow 2^{\R}$ be a maximal monotone graph, $a\in\R,\ \varepsilon >0$. Then
\begin{equation*}
(\gamma+a)^Y_\varepsilon (r)=\gamma^Y_\varepsilon(r-\varepsilon a)+a\quad\text{for every }r\in\R.
\end{equation*}
\end{lemm}
\begin{pr}
Recalling the definitions of Yosida approximation and of the inverse operator of a multivalued operator (see, e.g., \cite{bar,bre}), for $r\in\R$ we have that
\begin{equation}
\label{pc1}
(\gamma+a)^Y_{\varepsilon}(r)=\frac{1}{\varepsilon}(I-(I+\varepsilon(\gamma+a))^{-1})(r)=\frac{r}{\varepsilon}-\frac{y}{\varepsilon} ,
\end{equation}
\noindent where $y$ is such that $r\in y+\varepsilon(\gamma+a)(y)=y+\varepsilon\gamma(y)+\varepsilon a$. We infer
$$
r-\varepsilon a\in y+\varepsilon\gamma(y),
$$
\noindent which is equivalent to  $y=(I+\varepsilon\gamma)^{-1}(r-\varepsilon a)$.
On the other hand, we have
\begin{equation*}
\gamma^Y_\varepsilon(r-\varepsilon a)+a=\frac{1}{\varepsilon}(I-(I+\varepsilon\gamma)^{-1})(r-\varepsilon a)+a
=\frac{r}{\varepsilon}-a+\frac{1}{\varepsilon}y+a=\frac{r}{\varepsilon}-\frac{y}{\varepsilon} ,
\end{equation*}
whence the thesis follows from \eqref{pc1}.
\end{pr}

Now, in view of definitions \eqref{yosidaint} and \eqref{yosidabound}, we prove that a variation of  \eqref{condcomp} holds true at the $\varepsilon$-level for our Yosida regularizations. Then, we \rosso{will} take advantage of it in the following \rosso{estimates}.

\begin{lemm} \label{growcondeps}
There holds
\begin{equation} \label{eqgrowcondeps}
|\beta_\varepsilon (r)|\leq \eta|\beta_{\Gamma,\varepsilon}(r)|+C_\Gamma\quad\text{  for every }r\in\R.
\end{equation}
\end{lemm}
\begin{pr}
It is trivial to check \eqref{eqgrowcondeps} for $r=0$.
Then, we distinguish the cases $r>0$ and $r<0$.
We observe that \eqref{convfunct} and \eqref{domcont}
imply that $r,\ \beta_{\Gamma}^{\circ}(r)$ and $\beta^{\circ}(r)$
have the same sign for every $r\in D(\beta_{\Gamma})$. If $r> 0$ and $r\in D(\beta_{\Gamma})$, we can multiply inequality \eqref{condcomp} by $\varepsilon
\, ( >0) $ and then add to both sides $r$, thus obtaining the equivalent inequality
\begin{equation*}
r+\varepsilon\beta^{\circ}(r)\leq \varepsilon(\eta
\beta^{\circ}_\Gamma(r)+C_\Gamma)+r,
\end{equation*}
\noindent which, in terms of operators, can be rewritten as
$
(I+\varepsilon\beta^{\circ})(r)\leq (I+\varepsilon(\eta\beta_\Gamma^{\circ}+C_\Gamma))(r).
$
Hence, thanks to the maximal monotonicity of the graphs $\beta$ and $\eta\beta_\Gamma+C_\Gamma$, we infer that for their resolvent operators the following inequality holds:
\begin{equation} \label{in3}
(I+\varepsilon\beta)^{-1}(s)\geq (I+\varepsilon(\eta\beta_\Gamma+C_\Gamma))^{-1}(s)\quad\text{for every }s> \varepsilon C_\Gamma.
\end{equation}
Recalling the definition of Yosida regularization of a
maximal monotone graph given in \eqref{yosida}, we observe that \eqref{in3} reduces to
\begin{equation} \label{in4}
(I-\varepsilon\beta^Y_\varepsilon)(s)\geq (I-\varepsilon(\eta\beta_\Gamma+C_\Gamma)^Y_\varepsilon)(s),
\end{equation}
whence we obtain
$$
\beta^Y_\varepsilon(s)\leq (\eta\beta_\Gamma+C_\Gamma)^Y_\varepsilon(s).
$$
Applying Lemma~\rosso{\ref{yosidasum}} and Lemma~\rosso{\ref{yosidacost}} to the right-hand side of the above inequality, we deduce that
$$
(\eta\beta_\Gamma+C_\Gamma)^Y_\varepsilon(s)=(\eta\beta_\Gamma)^Y_\varepsilon(s-\varepsilon C_\Gamma)+C_\Gamma=\eta(\beta_\Gamma)^Y_{\varepsilon\eta}(s-\varepsilon C_\Gamma)+C_\Gamma
$$
and consequently
\begin{equation*}
\beta^Y_\varepsilon(s)\leq\eta(\beta_\Gamma)^Y_{\varepsilon\eta}(s-\varepsilon C_\Gamma)+C_\Gamma\quad\text{for every }s> \varepsilon C_\Gamma .
\end{equation*}
Hence, we infer
\begin{equation} \label{in6}
\beta^Y_\varepsilon(s+\varepsilon C_\Gamma)\leq\eta(\beta_\Gamma)^Y_{\varepsilon\eta}(s)+C_\Gamma\quad\text{for every }s> 0.
\end{equation}
At this point, we consider the left-hand side of \eqref{in6} and note that,
by the monotonicity of $\beta^Y_\varepsilon$,
$
\beta^Y_\varepsilon(s)\leq\beta^Y_\varepsilon(s+\varepsilon C_\Gamma).$ Then,
owing to \eqref{yosidaint}-\eqref{yosidabound}, we have that
$$
\beta_\varepsilon(s)\leq\eta\beta_{\Gamma,\varepsilon}(s)+C_\Gamma
$$
\noindent for every $s>0$. Arguing similarly for $r<0$, we finally show \eqref{eqgrowcondeps}.
\end{pr}

\subsubsection{Third a priori estimate}

We test now equation \eqref{varformeps} by $\beta_\varepsilon(u_\varepsilon(t))$
with $t\in(0,T)$, then we integrate over $(0,T)$. \rosso{Letting $\delta$ be} an arbitrary positive number and adding the positive term $\delta\int_\Sigma|
\beta_\varepsilon(u_\varepsilon)|^2$ to both sides for convenience, we rearrange
and obtain
\begin{align}
&  \int_{\Omega}\hat{\beta}_\varepsilon(u_\varepsilon (T)) +\int_{\Gamma}\hat{\beta}_\varepsilon(u_\varepsilon(T))+\int_Q \beta'_\varepsilon(u_\varepsilon)|\nabla u_\varepsilon|^2
\label{firstest}\\
& +\nu\int_{\Sigma}\beta'_\varepsilon(u_\varepsilon)|\nabla_\Gamma u_\varepsilon|^2+\int_Q |\beta_\varepsilon(u_\varepsilon)|^2+\delta\int_\Sigma|\beta_\varepsilon(u_\varepsilon)|^2 \notag\\
& =\int_\Omega\hat{\beta}_\varepsilon(u_0)+\int_\Gamma\hat{\beta}_\varepsilon(u_0)+\int_Q(f-\pi(u_\varepsilon))\beta_\varepsilon(u_\varepsilon)\notag\\
& +\int_\Sigma (f_\Gamma-\pi_\Gamma(u_\varepsilon))\beta_\varepsilon(u_\varepsilon)+\int_\Sigma\left(\delta|\beta_\varepsilon(u_\varepsilon)|^2- \beta_\varepsilon(u_\varepsilon)\beta_{\Gamma,\varepsilon}(u_\varepsilon)\right). \notag
\end{align}
The first two terms on the right-hand side of \eqref{firstest} can be easily treated owing to \eqref{ineqeps} and \eqref{betacondex}, thus obtaining
\begin{equation*}
\int_\Omega \hat{\beta}_\varepsilon(u_0)+\int_\Gamma\hat{\beta}_\varepsilon(u_0) \leq C.
\end{equation*}
In view of \eqref{datacond}, \eqref{lipschitz} and
\eqref{firstaprepsn}, by Young's inequality we infer that
\begin{align*}
&\int_Q(f-\pi(u_\varepsilon))\beta_\varepsilon(u_\varepsilon)\leq C + \frac{1}{2}\norm{\beta_\varepsilon(u_\varepsilon)}^2_{L^2(0,T;H)}
\end{align*}
and
\begin{equation*}
\int_\Sigma (f_\Gamma-\pi_\Gamma(u_\varepsilon))\beta_\varepsilon(u_\varepsilon) \leq \frac{\delta}{2}\norm{\beta_\varepsilon(u_\varepsilon)}^2_{L^2(0,T;H_\Gamma)} + c_\delta.
\end{equation*}
We have now to handle the last integral on the right-hand side of \eqref{firstest}. To this aim, we observe that Lemma \ref{growcondeps} and Young's inequality entail
\begin{align*}
&\delta|\beta_\varepsilon(u_\varepsilon)|^2-\beta_\varepsilon(u_\varepsilon)\beta_{\Gamma,\varepsilon}(u_\varepsilon)=
\delta|\beta_\varepsilon(u_\varepsilon)|^2-|\beta_\varepsilon(u_\varepsilon)||\beta_{\Gamma,\varepsilon}(u_\varepsilon)| \\
& \leq |\beta_\varepsilon(u_\varepsilon)|^2(\delta-1/\eta)+C_\Gamma|\beta_\varepsilon(u_\varepsilon)|/\eta \\
& \leq |\beta_\varepsilon(u_\varepsilon)|^2(\delta-1/2\eta)+C^2_\Gamma/2\eta\leq C^2_\Gamma/2\eta\quad\text{ a.e. on }\Sigma
\end{align*}
\noindent whenever $\delta<1/2\eta$.
Thus, by fixing  $\delta<1/2\eta$ it turns out that
the last integral of \eqref{firstest} is uniformly bounded.  Then, \blu{neglecting} the first four positive terms on the left-hand side, from  \eqref{firstest} it follows that
\begin{equation} \label{firstestimate}
\norm{\beta_\varepsilon(u_\varepsilon)}_{L^2(0,T;H)} +
\norm{\beta_\varepsilon(u_\varepsilon)}_{L^2(0,T;H_\Gamma)}
\leq C.
\end{equation}
The first part of estimate \eqref{firstestimate} has some consequences on the regularity of $u_\varepsilon $ and of the boundary term $\partial_n u_\varepsilon$.
Indeed, taking an arbitrary $ z  \in H^1_0(\Omega) $ as test function in \eqref{varformeps}, we can recover the partial differential equation satisfied by $ u_\varepsilon $ in $Q$, which can be written as
\begin{equation} \label{ossreg}
-\Delta u_\varepsilon = f-\partial_t u_\varepsilon-\beta_\varepsilon(u_\varepsilon)-\pi(u_\varepsilon)
\end{equation}
\noindent and holds (for instance) in the sense of distributions on $Q$. Now,
thanks to the estimates \eqref{unibound} and \eqref{firstestimate} and in view
of the assumption \eqref{datacond} on $f$ and of
the Lipschitz continuity of $\pi$,
by comparison in \eqref{ossreg} we deduce that
\begin{equation} \label{laplacl2}
\norm{-\Delta u_\varepsilon}_{L^2(0,T;H)}\leq C .
\end{equation}
Hence, recalling \eqref{unibound} and using e.g. \cite[Theorem~3.1, p.~1.79]{brez}, applied in this case with
\begin{equation*}
A=-\Delta,\quad g_0=u_\varepsilon|_\Gamma,\quad p=2,\quad  r=0,\quad t=1,\quad s=\frac{3}{2} , \label{valpar}
\end{equation*}
we deduce the estimate
\begin{equation}
\label{pc2}
\int_0^T\norm{u_\varepsilon (t)}^2_{H^{3/2}(\Omega)}dt
\leq c \int_0^T\left(\norm{\Delta u_\varepsilon (t) }^2_{H}+\norm{
u_\varepsilon (t)  }^2_{V_\Gamma}\right)dt,
\end{equation}
which implies (cf.~\eqref{unibound})
\begin{equation}
\label{pc3}
\nu^{1/2} \norm{u_\varepsilon}_{L^2(0,T:H^{3/2}(\Omega))}\leq C .
\end{equation}
By virtue of \eqref{laplacl2}, \eqref{pc3} and
\cite[Theorem~2.27, p.~1.64]{brez} \blu{we are also able} to conclude that
\begin{equation}
\nu^{1/2} \norm{\partial_n u_\varepsilon}_{L^2(0,T; H_\Gamma) } \leq C.    \label{estnormder}
\end{equation}

\subsubsection{Fourth a priori estimate}
Using now \eqref{ossreg} and \eqref{estnormder}, we are able to recover from \eqref{varformeps} the variational formulation of the boundary equation
\begin{equation} \label{pc4}
\partial_n u_\varepsilon + \partial_t u_\varepsilon - \nu\Delta_\Gamma
u_\varepsilon + \beta_{\Gamma,\varepsilon}(u_\varepsilon) +
\pi_\Gamma(u_\varepsilon) = f_\Gamma
\end{equation}
on $\Sigma,$ that is
\begin{equation*}
\int_\Gamma \partial_n u_\varepsilon z+\int_\Gamma\partial_t u_\varepsilon z +\nu\int_\Gamma\nabla_\Gamma u_\varepsilon\cdot\nabla_\Gamma z+\int_\Gamma\beta_{\Gamma,\varepsilon}(u_\varepsilon)z+\int_\Gamma\pi_\Gamma(u_\varepsilon)z=\int_\Gamma f_\Gamma z
\end{equation*}
\noindent for every $z\in \mathcal{V}$, a.e. in $(0,T)$.  We aim to find an estimate for
$\beta_{\Gamma,\varepsilon}(u_\varepsilon)$. Then, we test the above equation by
$\beta_{\Gamma,\varepsilon}(u_\varepsilon)$ and integrate over $(0,T).$
We obtain
\begin{align}
& \int_\Gamma \hat{\beta}_{\Gamma,\varepsilon}(u_\varepsilon(T))
+\nu\int_\Sigma\beta'_{\Gamma,\varepsilon} (u_\varepsilon)|\nabla_\Gamma \,
u_\varepsilon|^2
+ \int_\Sigma|\beta_{\Gamma,\varepsilon}(u_\varepsilon)|^2  \label{thirdest} \\
& = \int_\Gamma \hat{\beta}_{\Gamma,\varepsilon}(u_0)  + \int_\Sigma \left(f_\Gamma-\pi_\Gamma(u_\varepsilon) - \partial_n u_\varepsilon\right)\beta_{\Gamma,\varepsilon}(u_\varepsilon). \notag
\end{align}
The first integral on the right-hand side of \eqref{thirdest} is bounded,
thanks to \eqref{ineqeps} and \eqref{betacondex}. In fact, we have
\begin{equation*}
\int_\Gamma \hat{\beta}_{\Gamma,\varepsilon}(u_0)\leq\big\Vert
\hat{\beta}_{\Gamma}(u_0)\big\Vert_{L^1(\Gamma)}.
\end{equation*}
For the last integral on the right-hand side we
invoke \eqref{datacond}, \eqref{lipschitz},
\eqref{firstaprepsn} and \eqref{estnormder}: then \blu{Young's} inequality
helps us to infer that
\begin{equation*}
\int_\Sigma \left( f_\Gamma-\pi_\Gamma(u_\varepsilon) -
\partial_n u_\varepsilon \right) \beta_{\Gamma,\varepsilon}(u_\varepsilon)
\leq \frac{1}{2} \int_\Sigma|\beta_{\Gamma,\varepsilon}(u_\varepsilon)|^2
+ C \left( 1+ \nu^{-1} \right).
\end{equation*}
Therefore, forgetting the two positive integrals on the left-hand side of \eqref{thirdest}, we plainly obtain
\begin{equation} \label{thirdestimate}
\nu^{1/2} \norm{\beta_{\Gamma,\varepsilon}(u_\varepsilon)}_{L^2(0,T;H_\Gamma)}\leq C.
\end{equation}

\subsection{Passage to the limit} \label{lim:it}

Thanks to the \blu{previous estimates}, we can now prove the existence of the solution, thus completing the proof of Theorem~\ref{exist}. Owing to \eqref{unibound}, \eqref{firstestimate} and \eqref{thirdestimate} and by standard compactness results, we infer that limit functions exist such that the following convergences hold true:
\begin{align}
 u_\varepsilon\stackrel{*}{\rightharpoonup} u\quad &\text{in
}L^\infty(0,T,V)\cap H^1(0,T;H) ,\label{conv1} \\
 u_\varepsilon|_\Gamma\stackrel{*}{\rightharpoonup} u|_\Gamma=v\quad &\text{in
}L^\infty(0,T;V_\Gamma)\cap H^1(0,T;H_\Gamma), \label{conv2}\\
\beta_\varepsilon(u_\varepsilon)\rightharpoonup\xi\quad &\text{in }L^2(0,T;H) ,
\label{conv3}\\
 \beta_{\Gamma,\varepsilon}(u_\varepsilon|_\Gamma)\rightharpoonup\xi_\Gamma\quad
&\text{in }L^2(0,T;H_\Gamma) \label{conv4}
\end{align}
\noindent at least for a not relabeled subsequence, as $\varepsilon\searrow 0$.

Now, we want to prove that the quadruplet $(u,v,\xi,\xi_\Gamma)$ is a solution
to our problem \eqref{reguniqu}-\eqref{varform}. Recalling that the embeddings
$V\subset H$ and $V_\Gamma\subset  H_\Gamma$ are compact, we can apply the
result in  \cite[Sect.~8, Cor.~4]{si} to infer that
\begin{align}
u_\varepsilon\rightarrow u\quad&\text{strongly in }C^0([0,T];H), \label{converend1} \\
u_\varepsilon|_\Gamma\rightarrow v\quad&\text{strongly in }C^0([0,T];H_\Gamma).
\label{converend2}
\end{align}
Note that the initial \blu{conditions \eqref{dataconduniq} are} fulfilled and moreover,
because of the Lipschitz continuity of $\pi$ and $\pi_\Gamma$, we have that
\begin{align*}
\pi(u_\varepsilon)\rightarrow\pi(u)\quad&\text{strongly in } C^0([0,T];H),\\
\pi_\Gamma(u_\varepsilon|_\Gamma)\rightarrow\pi_\Gamma(v)\quad&\text{strongly in }C^0([0,T];H_\Gamma) .
\end{align*}
Then, passing to the limit in \eqref{varformeps}, we recover the variational
equality \eqref{varform}. It remains to check \eqref{reguniqcsi}-\eqref{reguniqcsigamma}. Since \eqref{conv3} holds and
$ u_\varepsilon\rightarrow u $  strongly in $L^2(0,T;H) $, we easily deduce
$$
\limsup_{\varepsilon\searrow 0}\int_Q \beta_{\varepsilon}
(u_{\varepsilon})u_{\varepsilon}=\int_Q\xi u  ,
$$
\noindent which allows us to apply \cite[Prop.~1.1, Ch.~II]{bar} and
conclude that $\xi\in\beta(u)$ a.e. in $Q$. We can argue in a similar way for
$\beta_{\Gamma,\varepsilon}$ \blu{and} obtain  $\xi_\Gamma\in\beta_\Gamma(v)$
a.e. on $\Sigma$. Therefore, the proof of Theorem \ref{exist} is complete.

\begin{obs}[Further regularities] \blu{In our passage to the limit} we can
recover further regularities for the limit function $u$. In fact, from \eqref{laplacl2} we infer that $\Delta u$ belongs to $L^2(0,T;H)$ just
by the lower semicontinuity property of the norm.  Similarly, we deduce from
\eqref{estnormder} that $\partial_n u$ is in $L^2(0,T;H_\Gamma)$.
Now, since $L^2(0,T;H_\Gamma)\subset L^2(0,T;H^{-1/4}(\Gamma))$, we can
consider $\partial_n u$ as an element of $L^2(0,T;H^{-1/4}(\Gamma))$ and then \blu{read} the equation \eqref{boundan} as the elliptic equation
$$-\Delta_\Gamma v + v= F_\Gamma$$
where
$$
F_\Gamma : = \frac{1}{\nu}(f_{\Gamma}-\partial_n u-\partial_t v-\xi_\Gamma-\pi_{\Gamma}(v))+v \in L^2(0,T;H^{-1/4}(\Gamma))
$$
Using now the boundary version of \cite[Theorem~7.5, p.~204]{lions}, we can deduce that $v=u|_\Gamma$ belongs to $L^2(0,T;H^{2-1/4}(\Gamma))\subset L^2(0,T;H^{3/2}(\Gamma))$, whence also
$$
u\in L^2(0,T;H^2(\Omega))
$$
\noindent thanks to \cite[Theorem~3.1, p.~1.79]{brez}. Moreover, we deduce that
\begin{equation*}
\rosso{\partial_n u\in L^2(0,T;H^{1/2}(\Gamma)) .}
\end{equation*}
\end{obs}

\section{The case without the Laplace-Beltrami operator}  \label{case nu0}
\setcounter{equation}{0}

In this section we deal with the case in which the parameter $\nu$ of the problem, which so far has been assumed to be positive, is equal to $0$. Both the operators $\Delta_\Gamma$ and $\nabla_\Gamma$ formally disappear and then the corresponding contributions have to be ignored. In particular, the case $\nu=0$ is significant in one space dimension since the above boundary operators are meaningless in this case. Before stating our main results, we have to redefine the operator induced on $H_\Gamma$ by $\beta_\Gamma$, extending it in a suitable way.

\subsection{Extension of the subdifferential on the boundary}
We introduce now another boundary space
$$
W_{\Gamma}:=H^{1/2}(\Gamma)
$$
in place of $V_{\Gamma}$ and note that, consequently, the new space of the test functions (replacing $\mathcal{V}$) is nothing but $V$, thanks to the trace theorem. Recalling the identification $   H_\Gamma  \cong  H^{*}_\Gamma$, let \blu{us indicate by} $W^{*}_\Gamma$ the dual space of $W_\Gamma$, that is the space $W^{*}_\Gamma:=H^{-1/2}(\Gamma)$.

In view of the definition \eqref{convfunct} of the function $\hat{\beta}_\Gamma$, we introduce a suitable generalization of the induced operator $\beta_\Gamma$ acting on the boundary. To this aim, we associate to the function $\hat{\beta}_\Gamma$ \cmag{the functionals} $\hat{\beta}_{\Gamma_{H_\Gamma}}(\rosso{z})$ on $H_\Gamma$ and $\hat{\beta}_{\Gamma_{W_\Gamma}}(\rosso{z})$ on $W_\Gamma$ \blu{defined \cmag{by}}
\begin{align}
&\hat{\beta}_{\Gamma_{H_\Gamma}}(\rosso{z}):=\left\{
\begin{aligned}
&\textstyle \int_\Gamma \hat{\beta}_\Gamma(\rosso{z}) &\text{if }\rosso{z}\in H_\Gamma\text{ and }\hat{\beta}_\Gamma (\rosso{z}) \in L^1(\Gamma) \\
&+\infty &\text{if }\rosso{z}\in H_\Gamma\text{ and }\hat{\beta}_\Gamma (\rosso{z})\notin L^1(\Gamma)
\end{aligned}
\right. , \notag\\[0.3cm]
&\hat{\beta}_{\Gamma_{W_\Gamma}}(\rosso{z}):=\hat{\beta}_{\Gamma_{H_\Gamma}}(\rosso{z})\quad\text{if }\rosso{z}\in W_\Gamma. \notag
\end{align}
As it is well-known, $\hat{\beta}_{\Gamma_{H_\Gamma}} $ and $\hat{\beta}_{\Gamma_{W_\Gamma}} $ are convex and lower \rosso{semicontinuous} functionals on $H_\Gamma$ and $W_\Gamma$, respectively. They are also proper since $W_\Gamma$ contains all the constant functions. Now, we denote by $\beta_{\Gamma_{W^{*}_\Gamma}}:=\partial \hat{\beta}_{\Gamma_{W_\Gamma}}: W_\Gamma\rightarrow 2^{W^{*}_\Gamma}$ the subdifferential of $\hat{\beta}_{\Gamma_{W_\Gamma}}$, which is defined by
\begin{align}
 \rosso{\rho}\in \beta_{\Gamma_{W^{*}_\Gamma}}(\rosso{z})\
\text{ if and only if }
\ \rosso{\rho}\in W^{*}_\Gamma, \ \rosso{z}\in D(\hat{\beta}_{\Gamma_{W_\Gamma}})\ \text{ and }  \label{subdiffrel} \\
 \hat{\beta}_{\Gamma_{W_\Gamma}}(\rosso{z})\leq {_{W^{*}_\Gamma}}\langle\rosso{\rho},\rosso{z} -w\rangle_{W_\Gamma}+\hat{\beta}_{\Gamma_{W_\Gamma}}(w)\quad\forall w\in W_\Gamma.  \notag
\end{align}
This turns out to be a maximal monotone operator \blu{from $W_\Gamma$ to
\cmag{$W^{*}_\Gamma$}} satisfying $\beta_{\Gamma_{W^{*}_\Gamma}}(0)\ni 0$
(cf.~\eqref{convfunct}). We observe that such an operator is strictly related
to the operator $\beta_\Gamma$ employed in previous sections since (see,
e.g., \cite[Prop.~2.5]{bcgg})
$ \rosso{\rho}\in \beta_\Gamma(\rosso{z})$ in $H_\Gamma$ if and only if
$\rosso{\rho}\in\beta_{\Gamma_{W^{*}_\Gamma}}(\rosso{z})$ and $
\rosso{\rho}  \in H_\Gamma$. We emphasize that
in the framework $\nu=0$ we are no longer able to deal with the natural
extension to $H_\Gamma$ of $\beta_\Gamma$. The inclusion
\eqref{reguniqcsigamma} governing the dynamic of the phase is 
\blu{in fact reinterpreted by \eqref{subdiffrel}} in
the abstract setting of the duality pairing between $W^{*}_\Gamma$ and
$W_\Gamma$. Nevertheless, the physical consistence is somehow preserved since
such an inclusion forces the phase to
assume only meaningful values in the domain of
$\beta_{\Gamma_{W^{*}_\Gamma}}$; further, if
$\xi_\Gamma(t)\in\beta_{\Gamma_{W^{*}_\Gamma}}(v(t))$ for a.e. $t\in (0,T)$,
we have in particular that $v(t)\in D(\hat{\beta}_\Gamma)\text{ a.e. on
}\Gamma$.

\subsection{Properties of the solution}
Let us formulate the new problem. As we have done before, we ask the data to fulfill at first the general assumptions \eqref{datacond} and we look now for a quadruplet $(u,v,\xi,\xi_\Gamma)$ such that
\begin{align}
& u\in L^2(0,T;V)\cap H^1(0,T;V^{*}), \label{pc5}\\
& v:=u|_{\Gamma}\in L^2(0,T;W_\Gamma)\cap H^1(0,T;W^*_\Gamma) , \\
& \xi\in L^2(0,T;H)\ \text{ and }\ \xi\in\beta(u)\text{ a.e. in }Q ,\label{pc5bis}\\
& \xi_\Gamma\in L^2(0,T;W_\Gamma^{*})\ \text{ and }\ \xi_\Gamma(t)\in
\beta_{\Gamma_{W^{*}_\Gamma}}(v(t))\,\text{ for a.e.  }t\in(0,T),\label{pc5ter}\\
& u(0)=u_0,\quad v(0)=v_0  \label{pc5quater}
\end{align}
\noindent and satisfying for a.a. $t\in (0,T)$ and every $z\in V$
\begin{align}\label{varform0}
 \langle\partial_t u(t),z\rangle+\int_{\Omega}\nabla u(t)\cdot\nabla z +\int_{\Omega}(\xi(t)+\pi(u(t)))z+ {_{W^*_\Gamma}}\langle\partial_t v(t),z\rangle_{W_\Gamma}\\
{}+\int_{\Gamma}\pi_{\Gamma}(v(t))z+ {_{W^{*}_\Gamma}}\langle\xi_{\Gamma}(t),z\rangle_{W_\Gamma}=\int_{\Omega}f(t)z+\int_{\Gamma}f_{\Gamma}(t)z . \notag
\end{align}

We state now the modified result of continuous dependence on the data and uniqueness: we observe that this is very similar to the one given for the case $\nu>0$, the only difference being the disappearance of the term with the factor $\nu$. We \blu{do not present} the proof of the theorem because it follows faithfully the arguments developed in Section~\ref{con:dep}.

\begin{teor}[Continuous dependence and uniqueness]  \label{uniqnu0}
Under the assumptions \eqref{convfunct}-\eqref{subdiff},
\eqref{datacond}, \eqref{subdiffrel},
let $(u_1,v_1)$ and $(u_2,v_2)$ be the first two components of two different solutions of the problem \eqref{pc5}-\eqref{varform0}, each one corresponding to a set of data $u_{0,i},v_{0,i},f_i,f_{\Gamma,i}$ with $i=1,2$. Then, there is
a constant $c$ such that
\begin{align}
 &\norm{u_1-u_2}^2_{L^{\infty}(0,T;H)}
+\norm{\nabla(u_1-u_2)}^2_{L^2(Q)}
+\norm{v_1-v_2}^2_{L^{\infty}(0,T;H_{\Gamma})} \notag\\
& {}\leq c \, \Big\{\norm{u_{0,1}-u_{0,2}}^2_H + \norm{v_{0,1}-v_{0,2}}^2_{H_\Gamma}
 +\norm{f_1-f_2}^2_{L^2(0,T;H)}+\norm{f_{\Gamma,1}-f_{\Gamma,2}}^2_{L^2(0,T;H_{\Gamma})}\Big\} .  \notag
\end{align}
\noindent Moreover, any two solutions of the problem \eqref{pc5}-\eqref{varform0} with the same data necessarily coincide.
\end{teor}

\blu{In the following} we denote by $(u^{\nu}, v^{\nu},\xi^\nu,\xi^\nu_\Gamma)$ the solution of the problem in the case $\nu>0$ \cmag{(provided by Theorems~\ref{exist} and \ref{uniqueness})} and by $(u,v,\xi,\xi_\Gamma)$ the expected solution in the case $\nu=0$. \blu{Even in this case we have to reinforce our assumptions on the data requiring}:
\begin{align}
& u_0\in V,\quad \text{so that } \ v_0=u_0|_{\Gamma}\in W_{\Gamma} , \label{pc6}\\
& \hat{\beta}(u_0)\in L^1(\Omega)\ \ \text{and}\ \ \hat{\beta}_{\Gamma}(v_0)\in L^1(\Gamma).  \label{pc7}
\end{align}
Note that \eqref{pc6} is weaker than \eqref{initcondex} (cf.~\eqref{VH}).
Since we want to prove that the solution $(u,v,\xi,\xi_\Gamma)$ can be obtained as the asymptotic limit of the quadruplet $(u^{\nu}, v^{\nu},\xi^\nu,\xi^\nu_\Gamma)$, we have to approximate the initial data $u_0, \, v_0 $ in order to comply with the regularity in \eqref{initcondex}-\eqref{betacondex}.
Then, we introduce a sequence $\{ u_{0\nu} \}  $ such that
\begin{align}
& u_{0\nu} \in V \ \text{ and } \ v_{0\nu}=u_{0\nu}|_{\Gamma}\in V_{\Gamma} , \label{pc11}\\
& \blu{u_{0\nu} \rightharpoonup u_0 \, \text{ in } \, V , \ \text{ whence } \
v_{0\nu} \rightharpoonup v_{0}  \, \text{ in } \, W_{\Gamma} \,}\text{ as }\nu\searrow 0 , \label{pc12}\\
& \nu^{1/2} \norm{v_{0\nu}}^2_{ V_{\Gamma}} + \int_\Omega \hat{\beta}(u_{0\nu}) + \int_\Gamma \hat{\beta}_{\Gamma}(v_{0\nu}) \leq C.  \label{pc12bis}
\end{align}
Now, the question for the reader is whether such a sequence $\{ u_{0\nu} \} $
exists. Of course,  if $v_0\in V_\Gamma$ then we can take $u_{0\nu} = u_0$ for all $\nu>0$. However,
the answer is positive in general.
\begin{lemm}
\label{exuzeronu}
Under the assumptions \eqref{pc6}-\eqref{pc7}, there is a sequence
$\{u_{0\nu}\} $ such that $u_{0\nu}$ and $v_{0\nu}=u_{0\nu}|_{\Gamma}$
satisfy \eqref{pc11}-\eqref{pc12bis}.
\end{lemm}
\begin{pr}
One can choose $u_{0\nu}$ as the solution of the following variational equality
\begin{equation}
\label{pc31}
\int_{\Omega} ( u_{0\nu} -u_0) z + \nu \int_{\Omega}\nabla u_{0\nu} \cdot\nabla z + \nu \int_{\Gamma} \rho_{0\nu} z  =0
 \quad \hbox{for all } \, z\in V,
\end{equation}
where
\begin{equation}
\label{pc31bis}
\rho_{0\nu} \in \beta_\Gamma (v_{0\nu} )
 \quad \hbox{a.e. in } \, \Gamma  .
\end{equation}
Indeed, consider the function
\begin{equation}
\Phi (z) =
\begin{cases} \displaystyle
\frac12 \int_\Omega |\nabla z|^2 + \int_\Gamma \hat\beta_\Gamma
(z|_\Gamma )  & \text{  if $z\in V$ and $\hat\beta_\Gamma
(z|_\Gamma) \in L^1(\Gamma)$}\\[0.4cm]
+\infty &\text{ otherwise}
\end{cases} \ ,
\notag
\end{equation}
which is proper, convex and lower semicontinuos in $H$. The subdifferential
of $\Phi$ is a maximal monotone operator from $H$ to $H$ (see \cite[Example~4, pp.~63-67]{bar}) with domain
$$ \{ z\in H^2(\Omega ) \ : \ - \partial_n z \in \beta_\Gamma (z|_\Gamma )
 \quad \hbox{a.e. in } \, \Gamma \}. $$
Then, it is straightforward to deduce the existence and uniqueness
of $ u_{0\nu}\in  H^2(\Omega )  $ solving \eqref{pc31}-\eqref{pc31bis}.
Moreover, $ u_{0\nu} $ turns out to be the
solution to the boundary value problem
\begin{align}
\label{pc32} & u_{0\nu}  - \nu \Delta  u_{0\nu}   = u_{0}  \quad \hbox{in } \, \Omega, \\
\label{pc33} &  \partial_n  u_{0\nu} + \beta_\Gamma (v_{0\nu} ) \ni 0  \quad \hbox{on } \, \Gamma.
\end{align}
Now, taking $z= u_{0\nu} -u_0 $ in \eqref{pc31}, we easily recover the estimate
\begin{equation}
\label{pc34}
\frac1\nu \norm{u_{0\nu} -u_0 }^2_{ H}
+ \Phi (u_{0\nu})
\leq
\Phi(u_{0}),
\end{equation}
whence \eqref{pc11}-\eqref{pc12} and part of \eqref{pc12bis}, namely
$$
\int_\Gamma \hat{\beta}_\Gamma(v_{0\nu} )
\leq
C,
$$
follow. To complete the verification of \eqref{pc12bis}, we can take $z= \beta_\varepsilon (u_{0\nu})$ in \eqref{pc31}, recalling that $\beta_\varepsilon$ denotes the Yosida regularization of $\beta.$ We obtain
$$
\int_\Omega \hat{\beta}_\varepsilon (u_{0\nu} )
+ \nu  \int_\Omega {\beta}'_\varepsilon (u_{0\nu} ) |\nabla u_{0\nu}|^2
+ \nu \int_\Gamma \rho_{0\nu}  {\beta}_\varepsilon(v_{0\nu} )
\leq
\int_\Omega \hat{\beta}_\varepsilon (u_{0})
$$
and, observing that the second and third terms on the left-hand side are
nonnegative (see \eqref{pc31bis} and note that $0 \in \beta_\Gamma (0)$
and ${\beta}_\varepsilon (0)=0 $), by \eqref{ineqeps} \cmag{and \cite[Prop.~2.11, p.~39]{bre}} we have
$$
\int_\Omega \hat{\beta}_\varepsilon (u_{0\nu} )
\leq \int_\Omega \hat{\beta} (u_{0}) \quad \hbox{for all } \, \varepsilon \in (0,1), \,
\hbox{ so that } \,
\int_\Omega \hat{\beta} (u_{0\nu} ) \leq C.
$$
Finally, using the estimate (2.23) in \cite[p.~64]{bar}, we find out that
$$ \norm{  u_{0\nu} }_{H^2(\Omega)} \leq C \left(1+ \norm{ u_{0\nu}}_H + \nu^{-1}  \norm{u_{0\nu} -u_0 }_{ H} \right) \leq C \left(1+ \nu^{-1/2}  \right)$$
thanks to \eqref{pc34}. Then, it turns out that $\nu^{1/2} \norm{  v_{0\nu} }_{H^{3/2}(\Gamma)}  \leq C $ and property \eqref{pc12bis} is ensured.
\end{pr}

Next theorem shows the existence of \blu{the solution to the problem \eqref{pc5}-\eqref{varform0}}
along with suitable convergences of the components
$u^{\nu}, v^{\nu},\xi^\nu,\xi^\nu_\Gamma$ of the solution
to the problem with $\nu>0$ to the respective components $u,v,\xi,\xi_\Gamma$.
\begin{teor} [Existence and convergence as $\nu\searrow 0$] \label{exist0}
Assume  \eqref{convfunct}-\eqref{subdiff},  \cmag{\eqref{domcont}-\eqref{condcomp},}  \eqref{datacond},  \eqref{subdiffrel} \rosso{and \eqref{pc6}-\eqref{pc7}.}
Then there exists a quadruplet $(u,v, \xi, \xi_\Gamma)$ satisfying
\begin{align}
& u \in  L^{\infty}(0,T;V) \cap H^{1}(0,T;H), \label{existu0}\\
& v \in  L^{\infty}(0,T;W_\Gamma) \cap H^{1}(0,T;H_\Gamma) \label{existv0}
\end{align}
and solving problem~\eqref{pc5}-\eqref{varform0}. Moreover, denoting by
$(u^\nu,v^\nu,\xi^\nu,\xi^\nu_\Gamma)$ the solution to the
problem~\eqref{reguniqu}-\eqref{varform} \rosso{with initial data $u_{0\nu}, v_{0\nu}$ as in
\eqref{pc11}-\eqref{pc12bis},
the following convergences hold as $\nu\searrow 0$:}
\begin{align}
 u^\nu\stackrel{*}{\rightharpoonup} u\quad &\text{in }L^{\infty}(0,T;V) \cap
H^{1}(0,T;H) , \label{converu0}\\
 v^\nu\stackrel{*}{\rightharpoonup} v\quad &\text{in } L^{\infty}(0,T;W_\Gamma)
\cap H^{1}(0,T;H_\Gamma), \label{converv0}\\
 \xi^\nu \rightharpoonup \xi \quad &\text{in } L^2(0,T;H), \label{pc8}\\
\xi^\nu_\Gamma \rightharpoonup \xi_\Gamma\quad &\text{in }
L^2(0,T,W^*_\Gamma) . \label{pc9}
\end{align}
\end{teor}

\begin{pr}
We prove the convergences \eqref{converu0}-\eqref{pc9} by using compactness arguments. Let us recall the a priori estimates \eqref{firstaprepsn}, \eqref{unibound}, \eqref{firstestimate} and \eqref{laplacl2} \blu{holding} for the solution $(u^\nu, v^\nu, \xi^\nu,\xi^\nu_\Gamma)$. In particular, let us point out that
\begin{align} \label{pcue}
\norm{u^\nu}_{L^{\infty}(0,T;V) \cap H^1(0,T;H)} +\norm{v^\nu}_{  L^{\infty}(0,T;W_\Gamma)   \cap  H^1(0,T;H_\Gamma)} \qquad \qquad \\
+ \nu^{1/2} \norm{v^\nu}_{L^{\infty}(0,T;V_\Gamma)}
+ \norm{\xi^\nu }_{L^2(0,T;H)} \leq C. \notag
\end{align}
In view of \eqref{estnormder} and \eqref{thirdestimate},
it is clear that we cannot deduce a uniform bound for $\partial_n u^\nu$ and $ \xi^\nu_\Gamma$
in $L^2(0,T;H_\Gamma)$ independently of $\nu$.
On the other hand, with a similar reasoning and taking advantage of
\eqref{laplacl2}, we can apply once more \cite[Theorem~3.1, p.~1.79]{brez}
\blu{and} \cite[Theorem~2.27, p.~1.64]{brez}, where the new choice of
operators and parameters is \blu{now} given by
\begin{equation*}
A=-\Delta,\quad g_0=u^\nu|_\Gamma,\quad p=2, \quad r=0,\quad t=\frac{1}{2},\quad s=1.
\end{equation*}
Thus, we find the following estimate for the normal derivative
\begin{equation} \label{pcnd}
\norm{\partial_n u^\nu }_{L^2(0,T;{W^{*}_\Gamma})}\leq C .
\end{equation}
Next, let us consider the analog of \eqref{pc4}, i.e.,
\begin{equation} \label{pc10}
\partial_n u^\nu + \partial_t u^\nu - \nu\Delta_\Gamma
u^\nu + \xi_\Gamma^\nu +
\pi_\Gamma(u^\nu) = f_\Gamma  \quad \text{in } \, L^2(0,T;{H_\Gamma}),
\end{equation}
which holds in $W^{*}_\Gamma$, a.e. in $(0,T)$ as well.
Now, owing to \eqref{estnormder} and \eqref{thirdestimate},
multiplying  \eqref{pc10} by $\nu^{1/2}$, from a comparison of
terms we obtain
\begin{equation} \label{pc13}
\nu^{3/2} \norm{\Delta_\Gamma u^\nu}_{L^2(0,T;H_\Gamma)}\leq C.
\end{equation}
Due to \eqref{firstaprepsn} and to the fact that $\Delta_\Gamma$
is a linear and bounded operator from $V_\Gamma $ to $V_\Gamma^*
(\cong H^{-1} (\Gamma))$, we have
\begin{equation} \label{pc14}
\nu^{1/2} \norm{\Delta_\Gamma u^\nu}_{L^2(0,T;V^*_\Gamma)}\leq C.
\end{equation}
Then, combining \eqref{pc13} and  \eqref{pc14}, by interpolation (see, e.g.,
\cite[Theorem~2.20, p.~1.53]{brez}) we deduce that
\begin{equation} \label{pc15}
\nu \norm{\Delta_\Gamma u^\nu}_{L^2(0,T;W^*_\Gamma)}\leq C.
\end{equation}
Therefore, thanks to \eqref{pcnd} and \eqref{pc15},
a further comparison in \eqref{pc10} leads to
\begin{equation} \label{pc16}
\norm{\xi_\Gamma^\nu}_{L^2(0,T;W^{*}_\Gamma)}\leq C.
\end{equation}
Collecting then \eqref{pcue} and \eqref{pc16}, we conclude that there exists a
quadruplet $(u,v, \xi, \xi_\Gamma)$ such that the convergences
\eqref{converu0}-\eqref{pc9} and
\begin{equation} \label{pc17}
 \nu v^\nu\to 0 \quad \text{in }L^{\infty}(0,T;V_\Gamma)
\end{equation}
hold as $\nu \searrow 0$, in principle for a subsequence and
then for the whole family once we have checked that $(u,v, \xi, \xi_\Gamma)$
\blu{is} the unique solution to \blu{the} problem~\eqref{pc5}-\eqref{varform0}.

In order to verify that the limit quadruplet $(u,v, \xi, \xi_\Gamma)$ solves
problem~\eqref{pc5}-\eqref{varform0}, we can proceed as in
Subsection~\ref{lim:it}. In particular, we note that, as in
\eqref{converend1}-\eqref{converend2}, we still have
\begin{align}
u^\nu\rightarrow u\quad&\text{strongly in }C^0([0,T];H), \label{convmu1} \\
v^\nu \rightarrow v\quad&\text{strongly in }C^0([0,T];H_\Gamma),
\label{convmu2}
\end{align}
as the embedding $W_\Gamma \subset H_\Gamma$ is also compact.
Then, it is straightforward to obtain \eqref{varform0} when passing to the limit in
\eqref{varform}\rosso{; \blu{initial conditions \eqref{pc5quater} follow} from \eqref{convmu1}-\eqref{convmu2} and \eqref{pc12}}. Moreover, \eqref{pc5bis} can be deduced from,
e.g., \cite[Prop.~2.5, p.~27]{bre} thanks to \eqref{reguniqcsigamma},
\eqref{pc8}, \eqref{convmu1}, and
\begin{equation}\label{pc22}
\limsup_{\nu\searrow 0}\int_Q \xi^\nu
u^\nu \leq \int_Q \xi u  .
\end{equation}
We point out that \eqref{pc22} implies not only \eqref{pc5bis}
but the additional property
\begin{equation}\label{pc23}
\lim_{\nu\searrow 0} \, \int_Q \xi^\nu
u^\nu = \int_Q \xi u  .
\end{equation}
Analogously, one should prove that
\begin{equation}\label{pc24}
\limsup_{\nu\searrow 0}
\int_0^T {_{W^{*}_\Gamma}}\langle\xi^\nu_{\Gamma}(t),
v^\nu(t) \rangle_{W_\Gamma} dt
\leq
\int_0^T {_{W^{*}_\Gamma}}\langle\xi_{\Gamma}(t), v(t) \rangle_{W_\Gamma}
dt ,
\end{equation}
\rosso{in order} to infer \eqref{pc5ter}. \blu{From} \eqref{reguniqcsigamma} it follows
that \eqref{subdiffrel} holds for $\rosso{\rho}=\xi^\nu_{\Gamma}$ a.e. in $(0,T)$, whence
$$
\int_\Sigma \xi^\nu_\Gamma v^\nu = \int_0^T {_{W^{*}_\Gamma}}\langle\xi^\nu_{\Gamma}(t),
v^\nu(t) \rangle_{W_\Gamma} dt .
$$
Moreover, the convergences \eqref{pc9} and \eqref{converv0} are in force. Then, taking $z=u^\nu (t) $ in \eqref{varform} and integrating with respect to $t$, we obtain
\begin{align}
 \int_0^T {_{W^{*}_\Gamma}}\langle\xi^\nu_{\Gamma} (t),  v^\nu (t)\rangle_{W_\Gamma}dt
= \frac{1}{2}\norm{\rosso{u_{0\nu}}}^2_H -\frac{1}{2}\norm{u^\nu (T)}^2_H
- \norm{\nabla u^\nu }^2_{L^2(0,T;H)}
+ \int_Q(f-\pi(u^\nu ))u^\nu
\notag \\
{}-\int_Q\xi^\nu u^\nu +\frac{1}{2}\norm{\rosso{v_{0\nu}}}^2_{H_\Gamma} - \frac{1}{2}\norm{v^\nu (T)}^2_{H_\Gamma}
-\nu \norm{\nabla_\Gamma \, v^\nu }^2_{L^2(0,T;H_\Gamma)}
+\int_\Sigma(f_\Gamma-\pi_\Gamma(v^\nu )) v^\nu .
\notag
\end{align}
We pass now to the $\limsup $ as $\nu \searrow 0$, using \rosso{\eqref{pc12},} \eqref{convmu1}-\eqref{convmu2},  \eqref{converu0} \rosso{and}
the lower semicontinuity of the functional $z \mapsto \norm{\nabla z }^2_{L^2(0,T;H)}$ in the weak star topology of $ L^\infty(0,T;V) $,
the Lipschitz continuity of $\pi$ and $\pi_\Gamma$,  \eqref{pc23} and
the fact that
$$
-\nu \norm{\nabla_\Gamma \, v^\nu }^2_{L^2(0,T;H_\Gamma)} \leq 0 .
$$
Thus, it is straightforward to infer the inequality
\begin{align}
 \limsup_{\nu\searrow 0}
\int_0^T {_{W^{*}_\Gamma}}\langle\xi^\nu_{\Gamma} (t),  v^\nu (t)\rangle_{W_\Gamma}dt
\leq \frac{1}{2}\norm{u_0}^2_H -\frac{1}{2}\norm{u (T)}^2_H
- \norm{\nabla u }^2_{L^2(0,T;H)}
\notag \\
 \int_Q(f-\pi(u )-\xi )u
+ \frac{1}{2}\norm{v_0}^2_{H_\Gamma}
- \frac{1}{2}\norm{v (T)}^2_{H_\Gamma}
+\int_\Sigma(f_\Gamma-\pi_\Gamma(v )) v
\notag
\end{align}
and see that the right-hand side is nothing but $$\displaystyle \int_0^T
{_{W^{*}_\Gamma}}\langle\xi_{\Gamma}(t), v(t) \rangle_{W_\Gamma} dt ,$$
thanks to \eqref{varform0} \rosso{and \eqref{pc5quater}}. Hence \eqref{pc24} follows and Theorem~\ref{exist0} is completely proved.
\end{pr}

\section*{Acknowledgements}

The authors are deeply grateful to Antonio Segatti for the fruitful
discussions and \rosso{his encouragement during the preparation of
this paper. The financial support of
the MIUR-PRIN Grant 2008ZKHAHN  \emph{``Phase transitions, hysteresis
and multiscaling''} and of the IMATI of CNR in Pavia is gratefully
acknowledged.}

\end{document}